# ON THE BEHRENS–FISHER PROBLEM: A GLOBALLY CONVERGENT ALGORITHM AND A FINITE-SAMPLE STUDY OF THE WALD, LR AND LM TESTS

By Alexandre Belloni[1] and Gustavo Didier

*Duke University and Tulane University*

In this paper we provide a provably convergent algorithm for the multivariate Gaussian Maximum Likelihood version of the Behrens–Fisher Problem. Our work builds upon a formulation of the log-likelihood function proposed by Buot and Richards [5]. Instead of focusing on the first order optimality conditions, the algorithm aims directly for the maximization of the log-likelihood function itself to achieve a global solution. Convergence proof and complexity estimates are provided for the algorithm. Computational experiments illustrate the applicability of such methods to high-dimensional data. We also discuss how to extend the proposed methodology to a broader class of problems.

We establish a systematic algebraic relation between the Wald, Likelihood Ratio and Lagrangian Multiplier Test ($W \geq LR \geq LM$) in the context of the Behrens–Fisher Problem. Moreover, we use our algorithm to computationally investigate the finite-sample size and power of the Wald, Likelihood Ratio and Lagrange Multiplier Tests, which previously were only available through asymptotic results. The methods developed here are applicable to much higher dimensional settings than the ones available in the literature. This allows us to better capture the role of high dimensionality on the actual size and power of the tests for finite samples.

**1. Introduction.** The so-called Behrens–Fisher Problem may be straightforwardly stated as follows.

> Given two independent random samples $X_1, \ldots, X_{N_1}$ and $Y_1, \ldots, Y_{N_2}$, test whether their respective population means $\mu_1$ and $\mu_2$ coincide in the case where their covariances $\Sigma_1$ and $\Sigma_2$ are unknown.

Received April 2007; revised June 2007.
[1]Supported in part by the IBM Herman Goldstein Fellowship.
*AMS 2000 subject classification.* 62H15.
*Key words and phrases.* Behrens–Fisher Problem, high-dimensional data, hypothesis testing, algorithm, Wald Test, Likelihood Ratio Test, Lagrange Multiplier Test, size, power.

Despite the deceiving simplicity of its form, this problem has motivated a wealth of literature that began with the original works of Behrens [1] and Fisher [9, 10], and includes Welch [34, 35], Scheffé [26, 27], Yao [37], Robbins, Simons and Starr [23], Subrahmanian and Subrahmanian [33] and Cox [8], to name a few. For a review of the solutions for the BFP, see, for instance, Stuart and Ord [32] and Kim and Cohen [15]. The proposed solutions involve a myriad of different approaches, ranging from fiducial inference to Bayesian techniques.

In this paper we are interested in the classical multivariate version of the Behrens–Fisher Problem under Normality. In other words, $X_i$, $Y_j$ above should be interpreted as $d$-dimensional Gaussian random vectors with (vector) means $\mu_1$ and $\mu_2$, and $\Sigma_1$, $\Sigma_2$ as their respective $d \times d$ covariance matrices, where the sample sizes are greater than $d$. The sample covariance matrices are then positive definite (and thus invertible) with probability one if the true covariance matrices $\Sigma_1$ and $\Sigma_2$ are positive definite. Several applied problems can be formulated as Behrens–Fisher Problems (in particular, for high dimension) in diverse areas such as Speech Recognition (e.g., Chien [6]), Quality Control (e.g., Murphy [21]), Development Economics (e.g., Schramm, Renn and Biles [28]) and others.

In this context, the Likelihood Ratio Test is a natural choice in face of the well-known asymptotic behavior of the test statistic. It turns out, though, that the maximization of the log-likelihood function without restrictive assumptions on the covariances (e.g., $\Sigma_1 = \Sigma_2$) is a nontrivial matter. In general, explicit solutions to the maximization procedure do not exist, and due to nonconcavities in the objective function, the solution to the system of first order likelihood equations can lead to local optima, as shown in Buot and Richards [5]. Numerical algorithms are available in the literature (see, e.g., Mardia, Kent and Bibby [20] and Buot and Richards [5]), but their convergence properties are unknown.

The purpose of this paper is two-fold. First, to propose a provably convergent algorithm, called Cutting Lines Algorithm (CLA), for the Gaussian Maximum Likelihood Behrens–Fisher Problem (BFP, for short). Second, to use the algorithm to investigate the finite sample properties—size and power—of the Likelihood Ratio Test and of the asymptotically equivalent Wald and Lagrange Multiplier Tests in the context of the BFP. Such properties are generally unknown, especially in high-dimensional contexts.

The CLA avoids the trap of local maxima, which haunts most approaches in the literature, by aiming directly for the maximization of the log-likelihood function itself. For this purpose, we make use of the expression for the log-likelihood function recently proposed by Buot and Richards [5], which is particularly suitable for numerical methods.

The general maximization strategy may be schematically characterized as follows:

(i) Lift the log-likelihood maximization problem into a higher-dimensional setting by adding artificial variables and constraints. This new problem, the Lifted BFP, has the same solution as the original BFP;

(ii) Create a family of *convex* modifications (subproblems) of the Lifted BFP which we call Ellipsoidal Mean Estimation Problems (EMEP);

(iii) Solve a sequence of EMEP whose solutions (estimators of the mean) converge to the global solution of the Lifted BFP, that is, the proper maximum likelihood estimator of the mean.

Step (i) is a common procedure in Continuous Optimization when one wishes to find a simpler (but equivalent) description for the problem in a higher-dimensional setting.

Step (ii) generates a family of *convex* problems which is computationally tractable (in particular, first order conditions are not only necessary but also sufficient). In fact, due to the particular structure of the EMEP, we are able to propose a specialized method which solves each problem in this family very efficiently both theoretically and in (computational) practice.

Step (iii) plays the crucial role of avoiding local maxima to ensure the global optimality. To achieve that, the algorithm relies on the particular geometry of the nonconvexities associated with the problem. Such geometry allows for the construction of a sequence of approximations (based on supporting lines) to the log-likelihood function itself which can be efficiently optimized. We prove that the proposed method converges to a global solution. Furthermore, a simulation study provides strong numerical evidence of the suitability of the CLA for solving high-dimensional problems. Problems with dimension up to 1000 were solved in a couple of minutes.

We are particularly interested in the finite-sample properties of the Wald, Likelihood Ratio and Lagrange Multiplier Tests. We show that their respective test statistics satisfy systematic algebraic inequalities in the context of the BFP (such a result is known for classical linear models; see Savin [25], Berndt and Savin [2], and Breusch [4]). However, the CLA makes it possible to go one step further and provide a Monte Carlo study of the actual size and the power of such tests. Our results illustrate that the Wald Test is the most sensitive among the three to the impact of dimensionality, followed by the Likelihood Ratio Test. Especially when the sample size is (relatively) small with respect to the dimension, the Wald and the Likelihood Ratio Tests tend to over-reject the null hypothesis when we use the $\chi^2$ quantiles given by Wilks' Theorem. In contrast, the observed size of the Lagrange Multiplier Test seems to be rather robust with respect to dimensionality, with a slight tendency to under-reject the null hypothesis. Perhaps not surprisingly, these properties carry over to the power of the tests: for fixed sample sizes, the Wald Test displays higher power than the Likelihood Ratio Test, which in turn displays higher power than the Lagrange Multiplier Test. However, the similar shapes of the observed power curves of the three tests seem to suggest that, with appropriate test size adjustment, the three tests may end up showing similar power properties. We also applied the Bartlett correction to the Likelihood Ratio Test as proposed by Yanagihara and Yuan [36]. The corrected test tends to under-reject the null-hypothesis, especially for high-dimensional data. Accordingly, it usually displays lower power than the Lagrange Multiplier Test.

In recent years, interesting applied problems have been found which can be formulated, in a generic sense, in the framework of the Behrens–Fisher Problem under high dimension and low sample size (see, e.g., Srivastava [31]). However, in this case no tests invariant under nonsingular linear transformations exist (see Srivastava [30] and references therein). Thus, the classical Maximum Likelihood formulation of the Behrens–Fisher Problem does not seem appropriate. The case of high dimension and low sample size should probably be handled by different techniques (or through a nontrivial transformation to a new Behrens–Fisher Problem), and is a topic for future research.

The paper is organized as follows. Section 2 recasts the log-likelihood maximization problem as a nonconvex programming problem, and introduces the EMEP. Section 3 studies the geometry of the nonconvexities associated with the log-likelihood function. Section 4 presents the CLA and its convergence analysis. Section 5 studies the finite-sample properties of the Wald, Likelihood Ratio, Lagrange Multiplier and the Bartlett-corrected Likelihood Ratio Tests. It also contains a computational investigation of the properties of the CLA in comparison to some widely used heuristic methods. Section 6 conveys an extension of the analysis to general BFP-like problems. The Appendix contains the following: the pertinent Convex Analysis definitions; an explanation of the relation between the EMEP and the BFP; a special-purpose algorithm for solving the EMEP; and an alternative convergent algorithm, called Discretization Algorithm, for solving the BFP.

**2. Lifting and the EMEP.** Recall that our goal is to maximize the log-likelihood function of two independent random samples $\{X_i\}_{i=1}^{N_1}$ and $\{Y_i\}_{i=1}^{N_2}$, where $X_i \sim N(\mu, \Sigma_1)$ and $Y_j \sim N(\mu, \Sigma_2)$ are $d$-dimensional (random) vectors. From now on we assume that the sample covariance matrices $S_1$ and $S_2$ are invertible. The maximization problem means that we should find $\mu$, $\Sigma_1$ and $\Sigma_2$ that maximize

$$l(\mu, \Sigma_1, \Sigma_2) = -\frac{1}{2}\sum_{i=1}^{N_1}(X_i - \mu)'\Sigma_1^{-1}(X_i - \mu) - \frac{N_1}{2}\log\det\Sigma_1$$
(1)
$$-\frac{1}{2}\sum_{i=1}^{N_2}(Y_i - \mu)'\Sigma_2^{-1}(Y_i - \mu) - \frac{N_2}{2}\log\det\Sigma_2,$$

which is a highly nonlinear function of $\mu$, $\Sigma_1$ and $\Sigma_2$.

Recently, a more (computationally) tractable reformulation of (1) was proposed by Buot and Richards [5]. We restate it here as a lemma.

LEMMA 2.1 (Buot and Richards [5]). *Denote the vector sample means by*

$$\bar{X} = \frac{1}{N_1} \sum_{i=1}^{N_1} X_i \quad \text{and} \quad \bar{Y} = \frac{1}{N_2} \sum_{i=1}^{N_2} Y_i, \tag{2}$$

*and the sample covariance matrices by*

$$S_1 = \frac{1}{N_1} \sum_{i=1}^{N_1} (X_i - \bar{X})(X_i - \bar{X})' \quad \text{and} \quad S_2 = \frac{1}{N_2} \sum_{i=1}^{N_2} (Y_i - \bar{Y})(Y_i - \bar{Y})'. \tag{3}$$

*Assume $S_1$ and $S_2$ are invertible, and let $\widehat{\mu}$ be some possible value, or estimator, of $\mu$. The original problem of maximizing the likelihood function in $\mu$, $\Sigma_1$ and $\Sigma_2$ can be reduced to the minimization in $\widehat{\mu}$ of*

$$(1 + (\bar{X} - \widehat{\mu})' S_1^{-1} (\bar{X} - \widehat{\mu}))^{N_1/2} (1 + (\bar{Y} - \widehat{\mu})' S_2^{-1} (\bar{Y} - \widehat{\mu}))^{N_2/2}. \tag{4}$$

Expression (4) is already much more tractable than the original likelihood since it depends only on $\mu$. However, the likelihood maximization problem can become substantially more amenable to analysis if it is reformulated as a suitable mathematical programming problem. We can do that by *lifting* it to a higher-dimensional setting, that is, by including additional variables and constraints, and recasting it in the following way.

DEFINITION 2.1. The Lifted Gaussian Maximum Likelihood Behrens–Fisher Problem is to solve

$$\min_{\mu, u_1, u_2} f(u_1, u_2) = \frac{N_1}{2} \log(u_1) + \frac{N_2}{2} \log(u_2),$$
$$u_1 \geq 1 + (\bar{X} - \mu) S_1^{-1} (\bar{X} - \mu), \tag{5}$$
$$u_2 \geq 1 + (\bar{Y} - \mu) S_2^{-1} (\bar{Y} - \mu).$$

Since the solutions for the Lifted Gaussian Maximum Likelihood Behrens–Fisher Problem and the original Gaussian Maximum Likelihood Behrens–Fisher Problem must coincide, we will use the acronym BFP to refer to the former from now on.

The advantage to the lifting procedure is to confine the nonconvexity of the problem to just two variables, $u_1$ and $u_2$. Nevertheless, the objective function $f$ in (5) still poses a computational challenge since it is nonconvex. This means that we can still expect the existence of local solutions as suggested in [5], and further analysis is called for.

One may note, though, that $f$ is increasing in $u_1$ and $u_2$. Moreover, if one of the variables, say, $u_1$, is fixed, then the problem becomes fairly simple: for each value of $u_1$, we can obtain a solution $u_2^*(u_1)$. The same can be done with $u_1^*$ as a function of $u_2$. Therefore, associated with (5), we could think of a family of tractable "subproblems" (parameterized by $u_1$, e.g.). Next we will show how to relate the solutions to this family of subproblems to the solution of the original problem.

Let us focus on the constraints in (5). For a given $\widehat{\mu}$ (a "solution"), consider the squared Mahalanobis distance functions

(6)　$\mathcal{M}_{\bar{X}}(\widehat{\mu}) = (\bar{X} - \widehat{\mu})' S_1^{-1} (\bar{X} - \widehat{\mu})$　and　$\mathcal{M}_{\bar{Y}}(\widehat{\mu}) = (\bar{Y} - \widehat{\mu})' S_2^{-1} (\bar{Y} - \widehat{\mu})$.

Note the resemblance between such functions and the generalized distance function $G$ as defined in Kim [16]. They all give ellipsoids in $\widehat{\mu}$, but our use of the functions is different.

DEFINITION 2.2. The Ellipsoidal Mean Estimation Problem with respect to $X$ at level $v_1$ is to solve

(7)　　　　　　　　$h_X(v_1) := \min_{\mu} \{ \mathcal{M}_{\bar{Y}}(\mu) : \mathcal{M}_{\bar{X}}(\mu) \leq v_1 \}$

(analogously for $Y$).

In words, the EMEP with respect to $X$ at level $v_1$ is to find the *estimate* $\widehat{\mu}_{\text{EMEP}}$ of $\mu$ that minimizes the squared distance $\mathcal{M}_{\bar{Y}}$ under the constraint that the squared distance $\mathcal{M}_{\bar{X}}$ is bounded by $v_1$. The use of the word "estimate" can be justified in at least two ways. First, Gaussian maximum likelihood estimation is based upon finding a vector estimate $\widehat{\mu}_{\text{EMEP}}$ that minimizes a similar quadratic form. Second, the procedure above enjoys the reasonable property that if $\bar{X}$ and $\bar{Y}$ are close (in particular, equal), the solution $\widehat{\mu}_{\text{EMEP}}$ will also be close to $\bar{Y}$ (in particular, equal).

Even though the EMEP is simpler than the BFP, there is no closed-form solution for the former (for given $v_1$). Nonetheless, EMEP is, in fact, a *convex* problem and can be solved efficiently by a variety of available methods like gradient descent, interior-point methods, cutting-planes, and so on. Although all these methods are convergent and a few have good complexity properties (see [3, 13, 22]), in the Appendix we propose a specific algorithm which explores the particular structure of the problem. Not surprisingly, it enjoys better complexity guarantees and better practical performance than the aforementioned methods.

The BFP and the EMEP are, in fact, closely related. The BFP consists of achieving the optimal balance between the EMEP for $X$ and $Y$ *simultaneously*. This happens because the BFP is based upon the minimization of a function that is monotone in both distance functions. A precise characterization of the relation between the BFP and the EMEP is given in the following theorem.

THEOREM 2.1. *Let $(\widehat{\mu}, \widehat{u}_1, \widehat{u}_2)$ be a solution to the BFP. Then, $\widehat{\mu}$ is a solution to the EMEP with respect to $X$ (with respect to $Y$) at $\widehat{v}_1 = \mathcal{M}_{\bar{X}}(\widehat{\mu})$ [at $\widehat{v}_2 = \mathcal{M}_{\bar{Y}}(\widehat{\mu})$].*

PROOF. Without loss of generality, we will develop the argument only for the EMEP with respect to $X$.

Let $\widehat{\mu}_{\text{EMEP}}$ be a solution to the EMEP with respect to $X$ at some positive $v_1$. By the monotonicity of log, this means that the expression

$$\text{(8)} \qquad \frac{N_1}{2}\log(1+v_1) + \frac{N_2}{2}\log(1+\mathcal{M}_{\bar{Y}}(\mu))$$

is minimized at $\widehat{\mu}_{\text{EMEP}}$.

Now, let $(\widehat{\mu}, \widehat{u}_1, \widehat{u}_2)$ be a solution to the BFP problem. This means that the expression

$$\text{(9)} \qquad \frac{N_1}{2}\log(1+\mathcal{M}_{\bar{X}}(\mu)) + \frac{N_2}{2}\log(1+\mathcal{M}_{\bar{Y}}(\mu))$$

is minimized at $\widehat{\mu}$ and we have $\widehat{u}_1 = 1 + \mathcal{M}_{\bar{X}}(\widehat{\mu})$. Since expression (9) is an upper bound for expression (8) when we set $v_1 := \mathcal{M}_{\bar{X}}(\widehat{\mu})$, $\widehat{\mu}$ is also a solution to the EMEP with respect to $X$ at $v_1$. □

REMARK 2.1. Since $S_1$ and $S_2$ are positive definite matrices (not only semi-definite), for each level of $v_1$ the EMEP has a unique solution. However, this does not guarantee that the BFP also has a unique solution, since it could achieve the optimum at two different levels of the distance function.

**3. The underlying geometry of the lifted Behrens–Fisher Problem.** In this section we study the nature of the nonconvexities in (5), and we show how the feasible set is related to the EMEP. In particular, we obtain a convenient representation of the border of the feasible set that will be used in the algorithm developed in Section 4.

We start by considering the projection of the set of feasible points in (5) into the two-dimensional space of $u = (u_1, u_2)$:

$$\text{(10)} \qquad \mathcal{K} = \left\{ (u_1, u_2) \in \mathbb{R}^2 : \exists \mu \text{ such that } \begin{array}{l} u_1 \geq 1 + \mathcal{M}_{\bar{X}}(\mu) \\ u_2 \geq 1 + \mathcal{M}_{\bar{Y}}(\mu) \end{array} \right\}.$$

Figure 1 illustrates the geometry of $\mathcal{K}$. Since $\mathcal{M}$ is a convex function, $\mathcal{K}$ is a convex set. Also, $\mathcal{K}$ is unbounded, since $(u_1, u_2) \in \mathcal{K}$ implies that $(u_1 + \gamma_1, u_2 + \gamma_2) \in \mathcal{K}$ as well for arbitrarily values of $\gamma_1, \gamma_2 > 0$. Clearly, $u \in \mathcal{K}$ implies that $u_1 \geq 1$ and $u_2 \geq 1$.

Since the objective function of (5), $f(u) = f(u_1, u_2) = \frac{N_1}{2}\log(u_1) + \frac{N_2}{2}\log(u_2)$, depends only on the variables $u$, the optimal value of (5) equals

$$\text{(11)} \qquad \min\{f(u) : u \in \mathcal{K}\},$$

which still is a nonconvex minimization and potentially has many local minima.

However, the representation (11) has two desirable features. First, it completely separates the (nonconvex) minimization problem in two variables from the high dimensionality of $\mu$. This will be key to avoid the curse of dimensionality. Second, we can write out a compact region that contains the solution for (11). Define the following problem dependent constants:

$$
\begin{aligned}
\bar{L}_1 &= \min_{\mu}\{1 + \mathcal{M}_{\bar{X}}(\mu)\} = 1, \\
\bar{U}_2 &= \min_{u_2}\{u_2 : (\bar{L}_1, u_2) \in \mathcal{K}\} = 1 + \mathcal{M}_{\bar{Y}}(\bar{X}), \\
\bar{L}_2 &= \min_{\mu}\{1 + \mathcal{M}_{\bar{Y}}(\mu)\} = 1, \\
\bar{U}_1 &= \min_{u_1}\{u_1 : (u_1, \bar{L}_2) \in \mathcal{K}\} = 1 + \mathcal{M}_{\bar{X}}(\bar{Y}).
\end{aligned}
\tag{12}
$$

These quantities define a right triangle

$$
\{(\bar{L}_1, \bar{L}_2), (\bar{L}_1, \bar{U}_2), (\bar{U}_1, \bar{L}_2)\}, \tag{13}
$$

which contains the optimal solution $u^* = (u_1^*, u_2^*)$ for (11). In fact, observe that, by monotonicity, all points in $\mathcal{K}$ above or to the right of the hypotenuse of the triangle have a larger objective value than a point on the hypotenuse. Moreover, the remaining points of $\mathcal{K}$ are contained in the triangle. Therefore, the coordinates of the triangle vertices in (13) are lower and upper bounds on the optimal solution $(u_1^*, u_2^*)$, that is,

$$\bar{L}_1 \leq u_1^* \leq \bar{U}_1, \qquad \bar{L}_2 \leq u_2^* \leq \bar{U}_2.$$

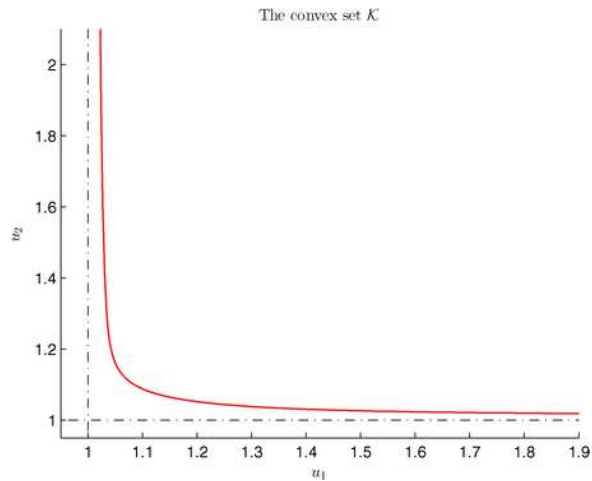

Fig. 1. *The convex set $\mathcal{K}$ consists of every point on and above the curve.*

In particular, if $\bar{X} = \bar{Y}$, the triangle degenerates into a single point (as pointed out in [5], the solution is trivial in this case).

Nevertheless, there is a representation cost associated with (11), in the sense that there is no closed-form representation for $\mathcal{K}$ involving only the variables $u$.

For this reason, we will make use of an additional function $g$ that gives information about (part of) the border of $\mathcal{K}$ (which is where the global optimum is expected to be found, given the quasi-concavity of $f$). The function $g$ is defined as

$$(14) \qquad g(u_1) := \min\{u_2 : (u_1, u_2) \in \mathcal{K}\}.$$

By construction, a point $(u_1, u_2)$ is in $\mathcal{K}$ if and only if $u_2 \geq g(u_1)$. It is easy to show that the function $g$ is convex (its epigraph is exactly the convex set $\mathcal{K}$) and decreasing in $u_1$.

Note that the function $g$ is directly related to the EMEP with respect to $X$ and the function $h_X$, since

$$(15) \qquad \begin{aligned} g(u_1) &= 1 + \min \mathcal{M}_{\bar{Y}}(\mu) = 1 + h_X(u_1 - 1), \\ u_1 - 1 &\geq \mathcal{M}_{\bar{X}}(\mu). \end{aligned}$$

In other words, evaluating $g$ at $u_1$ involves solving an EMEP with respect to $X$.

**4. An algorithm for the Behrens–Fisher Problem.** In this section we propose an algorithm, called the Cutting Lines Algorithm (CLA), that generates an $\varepsilon$-solution for the BFP. This means that the algorithm reports a feasible solution at which the objective function value lie within at most $\varepsilon$ from the value of the objective function at the optimal solution. Since the feasible solution is given for arbitrary $\varepsilon > 0$, convergence to an optimal solution holds.

The CLA builds upon a polyhedral approximation to the set $\mathcal{K}$. The method optimizes the objective function $f$ over $\widehat{\mathcal{K}}_k$ at each iteration. The minimizer point $(u_1, u_2) \in \hat{\mathcal{K}}_k$ is used to improve the polyhedral approximation for the next iteration.

As mentioned in the introduction, it is possible to propose an algorithm based upon the discretization of the range of values of $u_1$ where we need to evaluate $g(u_1)$. Such an algorithm, which we call a Discretization Algorithm (DA), can be proved to have better worst-case complexity guarantees than the ones obtained for the CLA. However, Section 5 shows that the practical performance of the CLA strongly dominates that of the DA, since the latter requires evaluating the function $g$—that is, solving an EMEP [see expression (15)]—at every point of the discretization. Thus, we focus on the CLA and defer the details of the DA to Appendix C.

4.1. *The cutting lines algorithm.* A good way to develop an algorithm for the BFP is to think of constructing sets that (i) approximate $\mathcal{K}$ and (ii) have a simple description involving $u$. Given the convexity of $\mathcal{K}$, polyhedral approximations to the set $\mathcal{K}$ are a natural candidate. Moreover, such approximations are rather convenient because it is simple to minimize the objective function $f$ over polyhedral sets in two dimensions (see Lemma 4.1 below).

4.1.1. *Building polyhedral approximations to $\mathcal{K}$.* Our sequence of polyhedral approximations will be based upon the function $g$. Given the results for the EMEP, relation (15) implies that, for any fixed value of $u_1$, not only can $g(u_1)$ be efficiently evaluated, but also a subgradient $s \in \partial g(u_1)$ (see Lemma B.1 for details) can be easily obtained. Suppose we choose a set of points $\{u_1^i\}_{i=1}^k$ and gather the triples

$$\{u_1^i, g(u_1^i), s^i\}, \qquad s^i \in \partial g(u_1^i), \qquad i=1,\ldots,k.$$

By the definition of subgradient, we have that

$$g(u_1) \geq g(u_1^i) + s^i(u_1 - u_1^i) \qquad \text{for all } i=1,\ldots,k \text{ and } u_1 \in \mathbb{R}.$$

Therefore, we can build a minorant polyhedral approximation $\widehat{g}_k$ for $g$ as follows:

(16) $$\widehat{g}_k(u_1) = \max_{1 \leq i \leq k} \{g(u_1^i) + s^i(u_1 - u_1^i)\}.$$

In turn, such a function can be used to build a polyhedral approximation for $\mathcal{K}$ defined as

$$\widehat{\mathcal{K}}_k = \{(u_1, u_2) \in \mathbb{R}^2 : u_2 \geq \widehat{g}_k(u_1)\}.$$

Figure 2 illustrates these relations.[1]

The advantage of working with the polyhedral approximation $\widehat{\mathcal{K}}_k$ instead of $\mathcal{K}$ is two-fold. First, $\widehat{\mathcal{K}}_k$ has a much nicer representation (via linear inequalities or extreme points) than $\mathcal{K}$ itself. This is particularly interesting for developing algorithms, which is our goal here. Second, as we anticipated, the minimization of the desired objective function $f(u_1, u_2) = \frac{N_1}{2}\log(u_1) + \frac{N_2}{2}\log(u_2)$ on $\widehat{\mathcal{K}}_k$ is rather tractable, as we show in the following lemma.

LEMMA 4.1. *Let $\widehat{\mathcal{K}}_k \subset \mathbb{R}^2_{++}$ be a (convex) polyhedral set. Then the function*

$$f(u_1, u_2) = \frac{N_1}{2}\log(u_1) + \frac{N_2}{2}\log(u_2)$$

*is minimized at an extreme point of $\widehat{\mathcal{K}}_k$.*

---

[1] Such approximation for convex sets can be traced back to the Cutting Planes Algorithm in the Optimization literature [3, 13, 14].

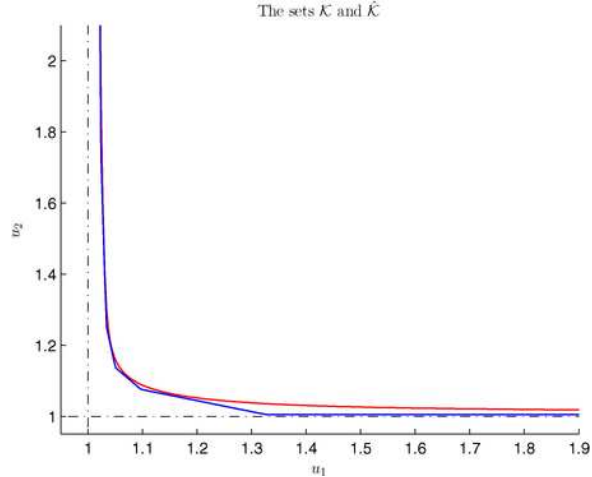

Fig. 2. *The convex set $\mathcal{K}$ and its outer polyhedral approximation $\widehat{\mathcal{K}}_k$. The extreme points of $\widehat{\mathcal{K}}_k$ are the kinks of the graph of the piecewise linear function $\hat{g}_k$.*

PROOF. First, note that since $\widehat{\mathcal{K}}_k \subset \mathbb{R}^2_+$, and because the nonnegative orthant is a pointed cone, $\widehat{\mathcal{K}}_k$ must have at least one extreme point. Second, the optimal solution cannot be an interior point of $\widehat{\mathcal{K}}_k$ (otherwise, we can strictly decrease both components simultaneously). Third, we recall that $f$ is a differentiable quasi-concave function. Therefore, its gradient is a supporting hyperplane for its upper level sets, which are convex.

Next, suppose that the minimum is achieved at a nonextreme point of $\widehat{\mathcal{K}}_k$, say, $x^* = \alpha z + (1-\alpha) y$, for $\alpha \in (0,1)$ and extreme points $z$, $y$. By the first order conditions, the gradient of $f$ induces a supporting line for $\mathcal{K}$ at $x^*$ on which both $z$ and $y$ lie. By the (strict) convexity of the upper level sets of $f$, $\min\{f(z), f(y)\} < f(x^*)$, a contradiction. $\square$

Since $\widehat{\mathcal{K}}_k$ is an outer approximation to $\mathcal{K}$, minimizing $f$ over $\widehat{\mathcal{K}}_k$ yields a lower bound on the optimal value of (5) for every $k$. Figure 3 illustrates the minorant approximation of $f(u_1, g(u_1))$ induced by $f(u_1, \hat{g}_k(u_1))$.

4.1.2. *The algorithm.* The CLA draws upon the minimization of the objective function over the polyhedral approximation $\widehat{\mathcal{K}}_k$ to $\mathcal{K}$, which, as shown in Lemma 4.1, needs to be carried out only over the extreme points of $\widehat{\mathcal{K}}_k$. A brief description of the algorithm follows. At iteration $k$, one has a set $f^i$, $i = 1, \ldots, k$, of values of the objective function at points $(u_1^i, u_2^i)$, $i = 1, \ldots, k$, respectively. The values $f^i$ are then compared to $\widehat{f}^k := f(\widehat{u}_1^k, \widehat{u}_2^k)$, where $(\widehat{u}_1^k, \widehat{u}_2^k)$ is the solution to the minimization of $f$ over $\widehat{\mathcal{K}}_k$. If the distance $\min_{0 \leq i \leq k}(f^i - \widehat{f}^k)$ is small enough (note that $f^i \geq \widehat{f}^k$), the algorithm stops.

Otherwise, it takes a new point $u_1^{k+1}$, slightly to the right of $\widehat{u}_1^k$, and generates its corresponding $u_2^{k+1} := g(u_1^{k+1})$ by solving an EMEP. The evaluation of the objective function $f$ at the pair $(u_1^{k+1}, u_2^{k+1})$ gives a new $f^{k+1}$, and the algorithm starts over.

| Cutting Lines Algorithm (CLA) | |
|---|---|
| **Input:** | Tolerance $\varepsilon > 0$, $u_1^1 = \min\{\bar{U}_1,\ (1+\varepsilon/N_1)\bar{L}_1\}$, $\widehat{g}_0 = 1$, $k = 1$. |
| **Step 1.** | Evaluate $u_2^k = g(u_1^k)$ and $s^k \in \partial g(u_1^k)$. Compute $f^k = \frac{N_1}{2}\log(u_1^k) + \frac{N_2}{2}\log(u_2^k)$. |
| **Step 2.** | Define $\widehat{g}_k(u_1) = \max_{0 \leq i \leq k}\{u_2^i + s^i(u_1 - u_1^i)\}$. |
| **Step 3.** | Compute $\widehat{f}_k = \min\{f(u_1, u_2): u_2 \geq \widehat{g}_k(u_1), u_1 \geq \bar{L}_1\}$ and the corresponding point $\widehat{u}^k = (\widehat{u}_1^k, \widehat{u}_2^k)$. |
| **Step 4.** | If $\min_{0 \leq i \leq k}(f^i - \widehat{f}^k) \leq \varepsilon$, report $\min_{0 \leq i \leq k} f^i$ and correspondent pair $(u_1^{i*}, u_2^{i*})$. |
| **Step 5.** | Else set $u_1^{k+1} \leftarrow \min\{\bar{U}_1,\ \widehat{u}_1^k(1+\varepsilon/N_1)\}$, $k \leftarrow k+1$, and goto Step 1. |

Note that each time a new iteration (say, $k+1$) starts, an updated polyhedral approximation $\widehat{\mathcal{K}}_{k+1}$ is constructed through the introduction of a new cut, based on the subgradient $\partial g(u_1^{k+1})$. A new cut removes one extreme point and creates at most two new extreme points. Therefore, the computational effort of minimizing $f$ over $\widehat{\mathcal{K}}_k$ grows only *linearly* with $k$ (in fact, by keeping track of previous evaluations, re-optimization can be done even faster).

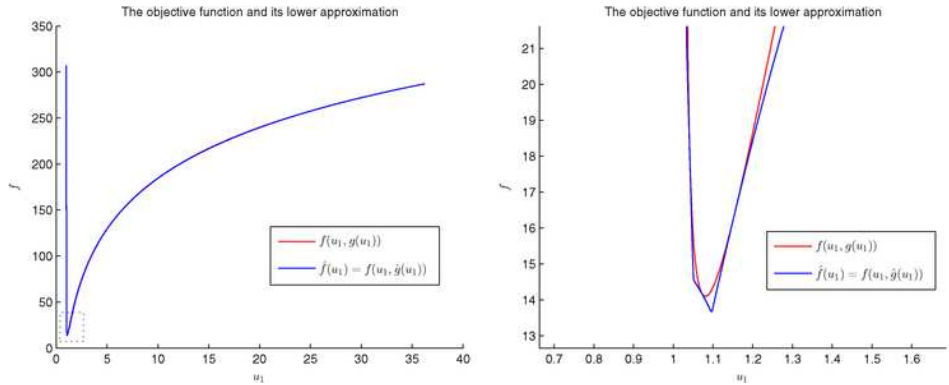

FIG. 3. *The outer polyhedral approximation for $\mathcal{K}$ leads to a minorant approximation for $f$. Therefore, lower bounds on the optimal value of (5) are derived if we minimize the minorant approximation $\widehat{f}$. The right figure is a zoom in on the dashed square area of the left figure.*

The next theorem shows that the CLA needs only a finite number of iterations to compute a $\varepsilon$-solution.

THEOREM 4.1. *The CLA reports an $\varepsilon$-solution to the original problem in at most $\lceil \frac{(\bar{U}_1 \bar{U}_2)(N_1 N_2)}{2\varepsilon^2} \rceil$ loops.*

PROOF. For $k \geq 1$, note that $u_1^{k+1} \leq \widehat{u}_1^k(1+\varepsilon/N_1)$, and suppose first that $u_2^{k+1} \leq \widehat{u}_2^k(1+\varepsilon/N_2)$. In this case, we have

$$f(u_1^{k+1}, u_2^{k+1}) = \frac{N_1}{2}\log(u_1^{k+1}) + \frac{N_2}{2}\log(u_2^{k+1})$$

$$\leq \varepsilon + \frac{N_1}{2}\log(\widehat{u}_1^k) + \frac{N_2}{2}\log(\widehat{u}_2^k)$$

$$= \varepsilon + \widehat{f}^k \leq \varepsilon + f^*,$$

and we have a $\varepsilon$-solution, since $(u_1^{k+1}, u_2^{k+1})$ is feasible.

Alternatively, if $u_2^{k+1} > \widehat{u}_2^k(1+\varepsilon/N_2)$, we have $u_2^{k+1} > 1$, which implies that $u_1^{k+1} < \bar{U}_1$. Therefore, $u_1^{k+1} = \widehat{u}_1^k(1+\varepsilon/N_1)$ and the next Cutting Lines approximation removes at least a rectangle of area $\frac{\varepsilon^2}{N_1 N_2}\widehat{u}_2^k \widehat{u}_1^k$ between the difference of $\widehat{\mathcal{K}}_k$ and $\mathcal{K}$. Since the area difference between these sets was bounded by $\bar{U}_1 \bar{U}_2/2$ at the very first iteration, the algorithm performs at most

$$\left\lceil \frac{(\bar{U}_1 \bar{U}_2)(N_1 N_2)}{2\varepsilon^2} \right\rceil$$

loops. □

This computational complexity result immediately yields the following convergence results.

COROLLARY 4.1. *For $\varepsilon_k \downarrow 0$, let $(u_1^k, u_2^k)$ be the $\varepsilon_k$-solutions to (11) and let the vectors $(\mu^k, u_1^k, u_2^k)$ be their induced $\varepsilon_k$-solutions to the (Lifted) BFP. Then, every accumulation point of the sequence $\{(\mu^k, u_1^k, u_2^k)\}_{k \in \mathbb{N}}$ is a solution to the BFP.*

COROLLARY 4.2. *The CLA can be used to generate a sequence of points that converge to a global solution to the (Lifted) Behrens–Fisher Problem.*

4.2. *Computational experiments with CLA and DA.* Our complexity bound for the CLA is worse than that for the DA. However, the DA solves the EMEP for every point of the discretized domain of $u_1$. In contrast, the CLA seeks to produce a certificate of $\varepsilon$-optimality at each iteration by comparing the best current solution and the solution to the minimization on $\widehat{\mathcal{K}}_k$.

TABLE 1
*Computational times (in seconds) and total number of iterations (which equal the number of EMEP problems solved) of the computational experiments with relative tolerance $\varepsilon = 10^{-3}$*

| Medium size instances | | | Average running times (in seconds) | | | Average iterations | |
|---|---|---|---|---|---|---|---|
| $d$ | $N_1$ | $N_2$ | Initialization | DA | CLA | DA | CLA |
| 20 | 100 | 200 | 0.01 | 4.45 | 0.006 | 6853.2 | 15.4 |
| 30 | 150 | 300 | 0.02 | 12.41 | 0.007 | 10859.6 | 17.5 |
| 40 | 200 | 400 | 0.03 | 12.92 | 0.009 | 9256.8 | 17.4 |
| 50 | 250 | 500 | 0.05 | 13.46 | 0.010 | 8414.6 | 18.5 |
| 60 | 300 | 600 | 0.08 | 23.71 | 0.012 | 10495 | 17.7 |
| 70 | 350 | 700 | 0.13 | 24.15 | 0.014 | 8502.1 | 17.6 |
| 80 | 400 | 800 | 0.18 | 42.40 | 0.020 | 9912.7 | 18.7 |
| 90 | 450 | 900 | 0.24 | 64.46 | 0.025 | 11796.4 | 18.3 |
| 100 | 500 | 1000 | 0.32 | 67.46 | 0.036 | 9859.5 | 19.0 |
| Large size instances | | | Average running times (in seconds) | | | Average iterations | |
| $d$ | $N_1$ | $N_2$ | Initialization | DA | CLA | DA | CLA |
| 200 | 1000 | 2000 | 2.07 | — | 0.23 | — | 20.8 |
| 300 | 1500 | 3000 | 6.64 | — | 0.66 | — | 19.8 |
| 400 | 2000 | 4000 | 16.08 | — | 1.66 | — | 20.1 |
| 500 | 2500 | 5000 | 43.35 | — | 3.13 | — | 21.2 |
| 600 | 3000 | 6000 | 56.62 | — | 5.71 | — | 21.5 |
| 700 | 3500 | 7000 | 87.88 | — | 6.88 | — | 22.0 |
| 800 | 4000 | 8000 | 142.05 | — | 12.71 | — | 20.9 |
| 900 | 4500 | 9000 | 455.23 | — | 23.59 | — | 22.1 |
| 1000 | 5000 | 10000 | 671.80 | — | 28.25 | — | 22.3 |

In computational practice, this drastically reduces the number of necessary iterations to find an $\varepsilon$-solution, as can be seen in Table 1 (this table was generated in the same way as the Monte Carlo study of the tests sizes, as described in Section 5.2 below). Each entry of running times and iterations in Table 1 is an average over ten instances.

Table 1 reflects the expected computational behavior of the methods. As the dimension increases, more effort is needed but the CLA is order of magnitudes faster than the DA, since the latter requires the complete discretization of the interval $[\overline{L}_1, \overline{U}_1]$. Such requirement of evaluating the function $g$ on $O(1/\varepsilon)$ different points (remember that the complexity analysis is exact in the case of the DA) seems to be a naive approach, indeed.

The polyhedral approximation used in the CLA provides a way of focusing the search on a promising region, a concept well exploited in the Optimization literature. Table 1 also illustrates the number of loops required by each algorithm in the test problems.

The number of loops performed by the Discretization Algorithm depends only on the precision $\varepsilon$, and on the problem dependent values of $\bar{L}_1$ and $\bar{U}_1$. On the other hand, these problem dependent quantities do not seem to affect the CLA. This points to the question of whether there exists a (better) complexity analysis for the CLA which might be independent of these quantities.

The implementation of the algorithms is a simple task in any programming package where matrix inversion and spectral decomposition subroutines for positive definite matrices are available (e.g., R, Matlab, etc.). The remaining algorithmic operations (binary search, computation of extreme points, stopping criterion, etc.) follow a relatively simple logic and do not involve potential numerical instabilities. We do not claim to have the most efficient implementation of the methods proposed here. Nevertheless, our numerical results show that the CLA is computationally efficient and scales quite nicely as the data dimension $d$ increases. The underlying reason is the certificate of optimality that the method is constructing on each iteration. The value $\widehat{f}_{\min}$ provides a lower bound for the optimal solution which is used to construct a stopping criterion. For a problem whose dimension is greater than one thousand, numerical approximations on the computation of the spectral decomposition are a potential limitation of the method to solve the EMEP proposed in Appendix B. An alternative approach is to compute an inverse matrix at each iteration of the EMEP, which will lead to a more robust implementation at the cost of additional running time (see [7] for details).

In our experiments we use medium and large size instances where the data dimension $d$ varies from 20 to 1000. The results were generated using a relative precision of $\varepsilon = 10^{-3}$. We report the average over ten different instances. The DA has proved to be too cumbersome for large instances.

**5. Finite sample properties of the Wald, Likelihood Ratio and Lagrange Multiplier tests through the CLA.** Three commonly used multivariate tests based upon the maximization of the log-likelihood function are the Wald $(W)$, Likelihood Ratio $(LR)$, and the Lagrange Multiplier $(LM)$ Tests. Define $\theta = (\mu_1, \mu_2, \Sigma_1, \Sigma_2)$. For a certain hypothesized restriction on the parameter space of means

$$H_0 : c(\mu_1, \mu_2) = q,$$

let $\widehat{\theta}$ denote the unrestricted MLE of $\theta$, and let $\widehat{\theta}_r$ denote the MLE under the restriction $H_0$, that is, the solution to the problem

$$\max_{\theta} \quad l(\theta)$$
$$\text{subject to} \quad c(\mu_1, \mu_2) = q.$$

The test statistics of interest are defined as

$$W = [c(\widehat{\mu}_1, \widehat{\mu}_2) - q]'(\text{Var}(c(\widehat{\mu}_1, \widehat{\mu}_2) - q))^{-1}[c(\widehat{\mu}_1, \widehat{\mu}_2) - q],$$
$$LR = -2(l(\widehat{\theta}_r) - l(\widehat{\theta}))$$

and

$$LM = e^T \widehat{G}_r [\widehat{G}_r^T \widehat{G}_r]^{-1} \widehat{G}_r^T e,$$

where

(17)
$$\widehat{G}_r = [\widehat{g}_{1,r}^x, \ldots, \widehat{g}_{N_1,r}^x, \widehat{g}_{1,r}^y, \ldots, \widehat{g}_{N_2,r}^y]^T,$$
$$\widehat{g}_{i,r}^x = \nabla_{\widehat{\theta}_r} \log f(x_i, \widehat{\theta}_r) \quad \text{and} \quad \widehat{g}_{i,r}^y = \nabla_{\widehat{\theta}_r} \log f(y_i, \widehat{\theta}_r),$$

($f$ is the multivariate density function in question) and $e$ is a vector of ones. In the context of the BFP, the restriction can be written as $\mu_1 - \mu_2 = 0$ and the $W$ test statistic has the explicit form

$$W = (\bar{X} - \bar{Y})'(S_1/N_1 + S_2/N_2)^{-1}(\bar{X} - \bar{Y}).$$

The $W$ Test—which is a pure significance test—bears the computational advantage of not requiring the solution to the problem of finding the restricted MLE estimator (however, see Section 5.2 below).

The $W$, $LR$ and $LM$ Tests are asymptotically equivalent under the null hypothesis. However, their behavior can be rather different in small samples, and their finite sample properties are usually unknown, except for a few particular cases (see, e.g., Greene [12] and Godfrey [11]). In this section we use the CLA to investigate and compare the finite sample properties—size and power—of these tests. In particular, we are interested in the sensitivity of the tests to dimensionality.

We emphasize that the CLA allows for the study of the properties of the tests in high-dimensional contexts. In contrast, the literature on the BFP typically overlooks the issue and reports results for small dimensional problems, typically smaller than $d = 6$ and in general no greater than $d = 10$.

5.1. *Conflict among criteria.* It is well known that the $W$, $LR$ and $LM$ statistics for testing linear restrictions in the context of classical linear models satisfy the inequalities $W \geq LR \geq LM$ (see Savin [25], Berndt and Savin [2], Breusch [4] and Godfrey [11]). Before turning to simulations, we show that such inequalities also hold in the case of the BFP.

THEOREM 5.1. *For the BFP,*

(18)
$$W \geq LR \geq LM.$$

PROOF. To show the first inequality, note that, using since $\log(1+\delta) \leq \delta$, we have
$$LR \leq c_0 = \min_{\mu} N_1(\bar{X} - \mu)S_1^{-1}(\bar{X} - \mu) + N_2(\bar{Y} - \mu)S_2^{-1}(\bar{Y} - \mu).$$
The optimal solution of the right-hand side is achieved at $\widehat{\mu}_0 = (N_1 S_1^{-1} + N_2 S_2^{-1})^{-1}(N_1 S_1^{-1} \bar{X} + N_2 S_2^{-1} \bar{Y})$. Using $\widehat{\mu}_0$, and the matrix identities
$$(A+B)^{-1} = A^{-1} - A^{-1}(A^{-1} + B^{-1})^{-1}A^{-1} = A^{-1}(A^{-1} + B^{-1})^{-1}B^{-1},$$
we prove that $c_0 = (\bar{X} - \bar{Y})'(S_1/N_1 + S_2/N_2)^{-1}(\bar{X} - \bar{Y}) = W$.

Let $\hat{\mu}$ be a solution for the BFP. After simplifications, the $LM$ statistic can be written as
$$LM = N_1(\bar{X} - \hat{\mu})'\widehat{\Sigma}_1^{-1}(\bar{X} - \hat{\mu}) + N_2(\bar{Y} - \hat{\mu})'\widehat{\Sigma}_2^{-1}(\bar{Y} - \hat{\mu}).$$
Next note that
$$(\bar{X} - \hat{\mu})'\widehat{\Sigma}_1^{-1}(\bar{X} - \hat{\mu}) = (\bar{X} - \hat{\mu})'S_1^{-1}(\bar{X} - \hat{\mu}) - \frac{[(\bar{X} - \hat{\mu})'S_1^{-1}(\bar{X} - \hat{\mu})]^2}{1 + (\bar{X} - \hat{\mu})'S_1^{-1}(\bar{X} - \hat{\mu})}$$
by using a rank-one update formula[2] for $\widehat{\Sigma}_1^{-1}$. The result follows by considering the term for $Y$ as well and noting that $\log(1+\delta) \geq \delta - \frac{\delta^2}{1+\delta}$. □

5.2. *Monte Carlo study of the size of the test.* Inequalities (18) imply that the rejection rate of the $W$ Test is greater than or equal to that of the $LR$ Test, which in turn is greater than or equal to that of the $LM$ Test. A more accurate understanding of the extent to which this influences the size and the power of such tests can be obtained through simulations.

We performed a Monte Carlo study of the finite-sample properties of the $W$, $LR$ and $LM$ tests at sizes $\alpha = 0.01, 0.05, 0.10$. The rejection regions were defined based upon Wilks' Theorem on the asymptotic $\chi_d^2$ distribution of the test statistic.

The study also includes the Likelihood Ratio statistic with the Bartlett correction
$$B := \left(1 - \frac{\widehat{c}_1}{N-2}\right)LR,$$
where
$$\widehat{c}_1 = \frac{\widehat{\psi}_1 - \widehat{\psi}_2}{d},$$
$$\widehat{\psi}_1 = \frac{N_2^2(N-2)}{N^2(N_1-1)}\{\operatorname{tr}(S_1\overline{S}^{-1})\}^2 + \frac{N_1^2(N-2)}{N^2(N_2-1)}\{\operatorname{tr}(S_2\overline{S}^{-1})\}^2,$$

---

[2] For invertible $M$ and a vector $v$, the inverse of the rank-one update of M by $vv'$ can be written as $(M + vv')^{-1} = M^{-1} - \frac{M^{-1}vv'M^{-1}}{1+v'M^{-1}v}$.

$$\widehat{\psi}_2 = \frac{N_2^2(N-2)}{N^2(N_1-1)} \{\text{tr}(S_1\overline{S}^{-1}S_1\overline{S}^{-1})\} + \frac{N_1^2(N-2)}{N^2(N_2-1)} \{\text{tr}(S_2\overline{S}^{-1}S_2\overline{S}^{-1})\},$$

and $\overline{S} = \frac{N_2}{N}S_1 + \frac{N_1}{N}S_2$.

The Bartlett correction as defined above provides an $O(N^{-2})$ approximation to the mean of the $\chi_d^2$ distribution (more details can be found in Yanagihara and Yuan [36]). We will refer to the $LR$ Test under the Bartlett correction as the $B$ Test.

To facilitate comparison with other works on the multivariate BFP (e.g., Yao [37], Subrahmaniam and Subrahmaniam [33], Kim [16] and Krishnamoorthy and Yu [18]), we performed tests for the low dimensional cases of $d = 2, 5$ and 10, but we also included the higher-dimensional cases of $d = 25, 50, 75$, 100 and 200. For each $d$, the sample sizes used were $N_1 = 5d, 10d, 20d$, and $N_2 = 2N_1$. For a given dimension size $d$, each covariance matrix $\Sigma_i$, $i = 1, 2$, was constructed by creating an initial matrix $M_i$ with $N(0,1)$ entries, and then setting $\Sigma_i = M_i M_i'$.

The results can be seen in Figure 4 (the actual numerical output can be found in Table 4 in the Appendix D). Each entry was generated using 10,000 runs. The $W$ and the $LR$ tests tend to over-reject the null hypothesis, while the $LM$ Test tends to slightly under-reject it. We kept constant the ratio between the number of observations and the dimension so that we can observe how the quality of the approximation behaves as the dimensionality of the problem grows. One may notice how sensitive the $W$ and the $LR$ Tests are to increases in the dimension. Only for the (relatively) large sample case $N_1 = 20d$ does the $LR$ Test have actual size fairly close to $\alpha$. On the other hand, the $W$ Test appears to demand even (relatively) larger samples. For instance, when $d = 100$ and $\alpha = 0.10$, even when $N_1 = 20d$, the $W$ test is off by 3.8 percentage points. The ease of computation of the $W$ test statistic appears to come at a considerable price in terms of the accuracy of the test.

In contrast with the $W$ and the $LR$ Tests, the $LM$ shows remarkable robustness with respect to dimensionality. For all $\alpha$, there does not appear to be any clear (say, monotonic) pattern of change on the actual test size with respect to increases in dimensionality, or maybe even sample size $N_1$.

For all values of $\alpha$ and different sample sizes, the $B$ Test is roughly as accurate as the $LM$ Test for low dimensional settings (roughly, $d \leq 20$). For $d > 20$, though, it grossly over-compensates the over-rejection rates of the $W$ Test, with the possible exception of the comparatively large sample sizes $N_1 = 20d$.

Figure 4 illustrates the above comments. Accordingly, the $W$ Test usually shows the steepest curve of dimension versus actual test size for different $N_1$, while the $LM$ Test displays approximately horizontal curves, especially for higher-dimensional settings.

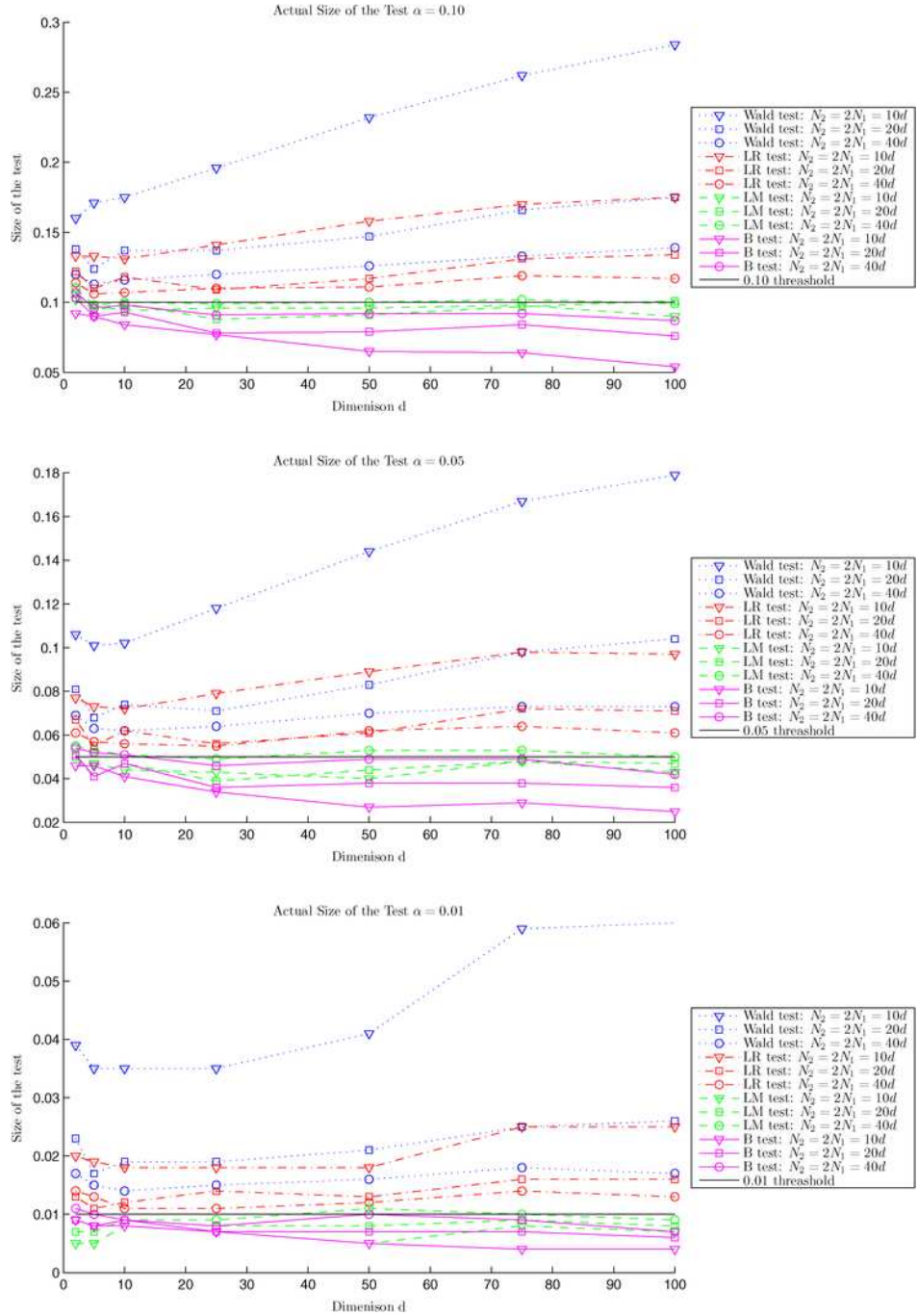

FIG. 4. *The behavior of the size of the tests when the dimension increases and the ratio between the number of observations and dimension is fixed.*

5.3. *Monte Carlo study of the power of the test.* We performed computational experiments on the power of the $W$, $LR$, $LM$ and $B$ Tests for the cases of dimension $d = 10, 50, 100$, and sample sizes $N_1 = 5d$, $10d$ and $20d$, with $N_2 = 2N_1$.

The analysis of the power for multivariate tests is naturally more difficult due to the multi-dimensionality of the parameter space. For this reason, we chose to investigate and compare the power of the $W$, $LR$, $LM$ and $B$ Tests over a standardized parameter space in the following sense. For each simulation run, covariance matrices $\Sigma_1$ and $\Sigma_2$ were (randomly) generated through the same procedure as the one for the evaluation of the sizes of the test. The mean of $X$, $\mu_1$, was set to zero by default. The choice of the mean(s) of $Y$, $\mu_2(\Delta)$, was made as solution(s) to the squared Mahalanobis distance equation(s)

$$(\mu_1 - \mu_2(\Delta))'(\Sigma_1 + \Sigma_2)^{-1}(\mu_1 - \mu_2(\Delta)) = \Delta^2,$$

where $\Delta$ represents a family of appropriately selected constants. For convenience, such solutions $\mu_2(\Delta)$ were always taken on some canonical axis, and the specific axis chosen changed across simulation runs. The use of randomly standardized Mahalonobis distances is justified by the fact that the BFP is defined without information on the population covariances.

The results are depicted in Figure 5, which contains plots for dimensions $d = 10$, 50 and 100. Colors represent tests, while geometric figures represent sample sizes (e.g., a triangle symbolizes $N_1 = 5d$).

Perhaps the most striking feature of all four plots ($d = 10$, 50 and 100) is the fact that, for a given sample size $N_1$, the *shapes* of the power curves for the four tests look alike. More specifically, given $N_1$, the curve for the $W$ Test looks like an up-shifted version of the curve for the $LR$ Test, which in turn looks like an up-shifted version of the curve for the $LM$ Test. The same is true for the curve for the $B$ Test, which lies mostly below the curve for the latter. The observed "order" of the curves should not come as a surprise. First, regarding the $W$, $LR$ and $LM$ Tests, because of the theoretical inequalities in Theorem 5.1. Second, because the simulation results for the test sizes show that the $W$ and $LR$ Tests tend to over-reject the null hypothesis (the former, substantially more than the latter), while the $LM$ Test has size close to $\alpha$ and the $B$ Test tends to under-reject the null hypothesis. In other words, we are essentially comparing tests of *different sizes* (see also the conclusions in Breusch [4] for the case of linear regression). The shape of the curves suggests the possibility that, if test size adjustment is made for the $W$ and $LR$ Tests, the power curves of the three tests may get rather close to each other. Such adjustment would imply, of course, going beyond Wilks' Theorem and developing exact quantiles, especially for the $W$ and the $LR$ Tests.

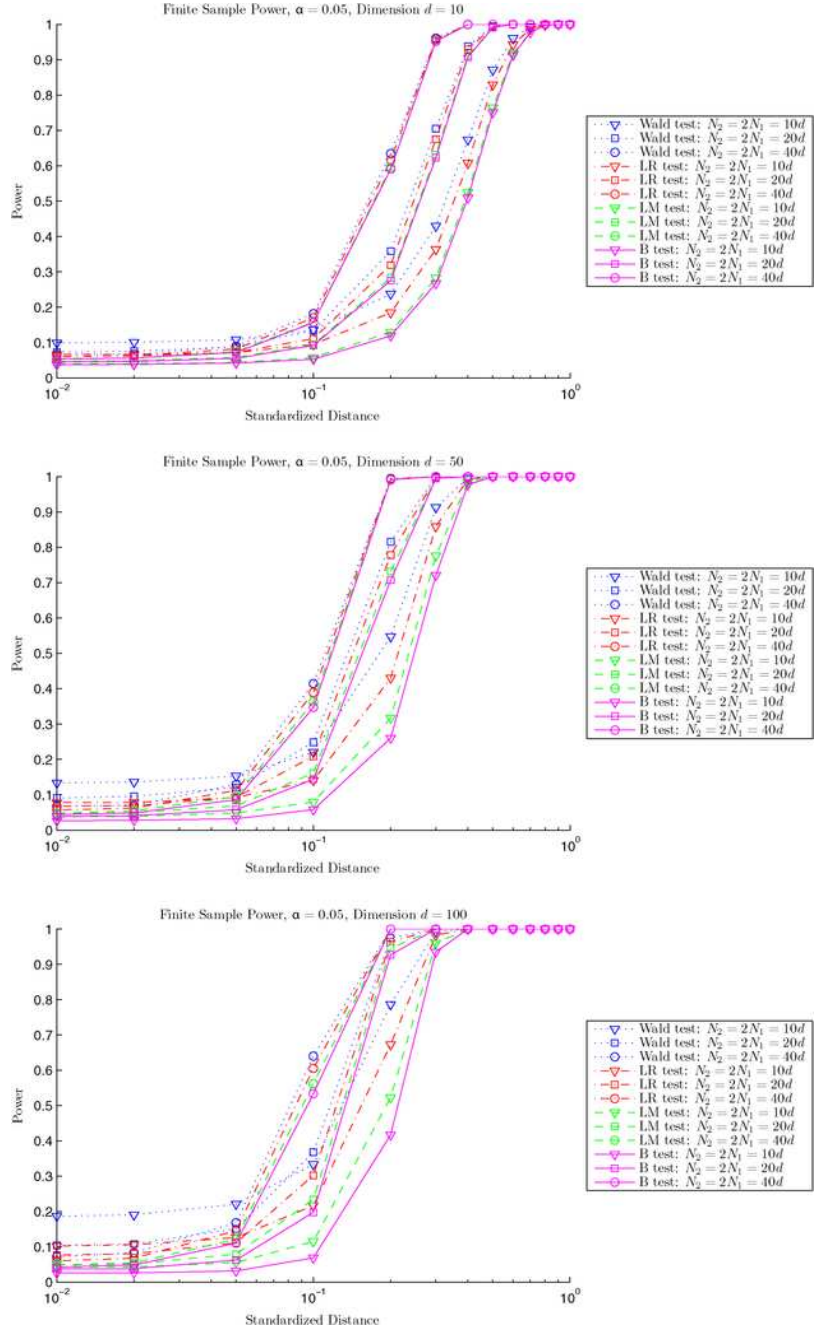

Fig. 5. *Monte Carlo study of the power of the W, LR, LM and B Tests for the size $\alpha = 0.05$ with different sample sizes and dimensions equal to 10, 50 and 100.*

The plot for the low-dimensional case of $d = 10$ displays a "well-behaved" pattern, in the sense that the curves for different tests and for the same sample size tend to be grouped together. In particular, the curves for sample size $N_1 = 40d$ are almost super-imposed, which means that, power-wise, the tests are nearly equivalent in this situation. Note that the curves for sample size $N_1 = 10d$ (triangle) lie above the remaining ones close to the origin, that is, in the case where the Mahalanobis distance between $\mu_1$ and $\mu_2$ is small. Again, this should not come as a surprise, since the simulation results for the test sizes (i.e., zero Mahalanobis distance between $\mu_1$ and $\mu_2$) show that relatively small sample sizes imply a tendency for over-rejection in the case of the $W$ and $LR$ Tests.

The effect of higher dimensionality can be seen in the two remaining plots ($d = 50$ and 100). The main impact seems to be greater vertical distances among the curves for the four tests, particularly for the cases of smaller sample sizes. Even for the higher-dimensional case $d = 100$, though, the larger sample size $N_1 = 20d$ brings the curves a lot closer to each other. As one might expect, larger sample sizes compensate for high dimension and point to the asymptotic equivalence of the $W$, $LR$, $LM$ and $B$ Tests.

5.4. *Performance of local methods/heuristics.* Up to the present, the numerical procedures applied to the multivariate Behrens–Fisher Problem have been heuristics or locally convergent methods. Since the CLA is a provably convergent method that constructs a certificate of global optimality, it provides a benchmark for the previous approaches. So, we are now able to address via Monte Carlo experiments the *statistically* important question of the performance of the $LR$ Tests based on some widely used heuristics vis-a-vis the $LR$ Test based on the CLA. Also, we are interested in the partially related issue the *computational* performance of these heuristics vis-a-vis the CLA.

There is a variety of different heuristics and it is usually hard (if not impossible) to make any general statement about them. However, in the case of the Behrens–Fisher Problem, we do have a "natural" initial point for these algorithms, that is, $\widehat{\mu}_0 := (N_1 S_1^{-1} + N_2 S_2^{-1})^{-1}(N_1 S_1^{-1} \bar{X} + N_2 S_2^{-1} \bar{Y})$. In fact, $\widehat{\mu}_0$ is at the same time: (i) from the algorithmic perspective, the solution to the first-order conditions of the objective function (the log-likelihood) with respect (only) to $\mu$, after we substitute $S_i$ for $\Sigma_i$, $i = 1, 2$; (ii) from the statistical perspective, the estimator of the mean $\mu$ associated with the $W$ statistic. Also, by the proof of Theorem 5.1, we have

$$W \geq LR_0 \geq LR, \tag{19}$$

where $LR_0$ is the log-likelihood ratio evaluated at $\widehat{\mu}_0$. Denote by $LR_h$ the (potentially suboptimal) Likelihood Ratio test statistic based upon a given heuristic method and with $\widehat{\mu}_0$ as its initial point. We can assume (through an *ad-hoc* modification of the heuristic, if necessary) that

$$(20) \qquad LR_0 \geq LR_h.$$

Thus, by (19), (20) and the fact that $LR_h \geq LR$, the statistic $LR_h$ also asymptotically follows a $\chi_d^2$ distribution. Nonetheless, the gap between the $W$ and the $LR$ statistics can be quite large in finite-samples (see Section 5).

Note that a $LR_h$ Test can only disagree with the $LR$ Test if the $LR_0$ Test rejects $H_0$ and the latter accepts $H_0$. We can perform Monte Carlo experiments (in the same way as in Section 5.2 for the tests sizes) to estimate how often the $LR_0$ and $LR$ Tests disagree. By (20), this provides a guarantee (i.e., an upper-bound) on the "discrepancy rate" between a $LR_h$ Test (i.e., based on any heuristic) and the $LR$ Test. The results in Table 2 show that this worst-case-scenario discrepancy rate is surprisingly small. The discrepancy rate for the $W$ Test with respect to the $LR$ Test—substantially higher—was also included for the sake of comparison.

Next we study the computational performance of three commonly used methods: Simulated Annealing (SA), Iterative Update (ItUp) and Newton's Method with Line Search (NM). (See, resp., [17, 19], [5] and [7] for discussions and implementation details of these methods.)

Table 3 reports the average performance of the heuristics with respect to running times, iterations and discrepancy rate based on their respective $LR_h$ Tests. As expected, the discrepancy rate is smaller than in Table 2, since the heuristics usually provide a solution superior to $\widehat{\mu}_0$. Not only that, the experiments suggest that both ItUp and NM are robust in terms of discrepancy rate (with $\widehat{\mu}_0$ as the initial point) even though they can be trapped in local minima. Moreover, their good (local) convergence properties are illustrated by the notably small number of iterations. On the other hand, SA seems to have trouble achieving local convergence, and its good performance with respect to errors appears to be a by-product of the chosen initial point. Not surprisingly, the convergence of the CLA turned out to be slower (i.e., larger number of iterations) than the local methods ItUp and NM. In fact, one should keep in mind that the CLA aims not only to find a good solution, but also to construct a *certificate* of global optimality, which is a much harder task. Regarding running time, the main computational cost of CLA is the spectral decomposition at initialization (see Table 1). The running time of the CLA after initialization is actually faster than SA, ItUp and NM at higher dimensions (cf. Tables 1 and 3).

TABLE 2
*Monte Carlo study of the discrepancy rates for $LR_0$ and $W$ Tests with respect to the $LR$ Test for different dimensions $d$ (for each entry, simulations were run until 50 "successes" were obtained)*

| $d$    | 2       | 5       | 10       | 20       | 30       | 40      |
|--------|---------|---------|----------|----------|----------|---------|
| $LR_0$ | 0.00176 | 0.00113 | 0.000932 | 0.000978 | 0.000825 | 0.00135 |
| $W$    | 0.0299  | 0.0325  | 0.0261   | 0.0468   | 0.0342   | 0.0513  |
| $d$    | 50      | 60      | 70       | 80       | 90       | 100     |
| $LR_0$ | 0.0009  | 0.0012  | 0.0012   | 0.0013   | 0.0017   | 0.0012  |
| $W$    | 0.0426  | 0.0702  | 0.0641   | 0.0796   | 0.0809   | 0.0935  |

TABLE 3
*The average performance of heuristics averaged over 5000 runs: time (seconds), iterations and discrepancy rate*

|     | **Algorithms** | | | | | | | | |
|-----|---|---|---|---|---|---|---|---|---|
|     | **SA** | | | **ItUp** | | | **NM** | | |
| $d$ | Time | Iter | Discrep | Time | Iter | Discrep | Time | Iter | Discrep |
| 2   | 0.02462 | 1000 | 0.001 | 0.00062 | 4.2 | 0 | 0.00118 | 2.3 | 0 |
| 5   | 0.02547 | 1000 | 0.001 | 0.00081 | 4.5 | 0 | 0.00135 | 2.6 | 0 |
| 10  | 0.02680 | 1000 | 0.001 | 0.00121 | 4.7 | 0 | 0.00166 | 2.6 | 0 |
| 20  | 0.03120 | 1000 | 0.001 | 0.00273 | 5.0 | 0 | 0.00310 | 2.7 | 0 |
| 30  | 0.04078 | 1000 | 0.001 | 0.00799 | 5.1 | 0 | 0.00662 | 2.9 | 0 |
| 40  | 0.04948 | 1000 | 0.001 | 0.00879 | 5.4 | 0 | 0.00787 | 3.1 | 0 |
| 50  | 0.06037 | 1000 | 0.001 | 0.01339 | 5.3 | 0 | 0.01179 | 3.0 | 0 |
| 60  | 0.07422 | 1000 | 0.001 | 0.01998 | 5.3 | 0 | 0.01726 | 3.0 | 0 |
| 70  | 0.09325 | 1000 | 0.001 | 0.03057 | 5.3 | 0 | 0.02537 | 3.0 | 0 |
| 80  | 0.11258 | 1000 | 0.003 | 0.04069 | 5.5 | 0 | 0.03408 | 3.1 | 0 |
| 90  | 0.13279 | 1000 | 0.001 | 0.05461 | 5.5 | 0 | 0.04414 | 3.0 | 0 |
| 100 | 0.16051 | 1000 | 0.001 | 0.07260 | 5.9 | 0 | 0.05852 | 3.5 | 0 |

We now make a few quick remarks regarding the implementation of the methods. First, all methods do require a matrix inversion routine: SA and CLA, only on the first iteration; ItUp and NW, on every iteration. Second, in contrast to CLA, ItUp and NW, the calibration of additional parameters is needed for the SA. Third, the implementation of SA and ItUp is very simple, while the Line Search for NM is slightly more difficult. Fourth, unlike the other methods, the CLA involves the additional implementation costs associated with the optimality certificate based on $\widehat{\mathcal{K}}_k$, and with the spectral decomposition of a positive definite matrix (see also Section 4.2).

**6. Extension to Behrens–Fisher-like Problems.** It should be noted that the methodology proposed in this paper can be applied to a much broader class of problems. Strictly speaking, all we need is to be able to replicate the strategy of constructing lifted problems whose solution lie on extreme points of a two dimensional convex domain,[3] and to evaluate the subproblems which define the convex domain. A sufficient condition for this is the quasi-concavity of the objective function of the lifted problem and the convexity of the subproblems.

To set up a broader framework, assume we have two random samples $\{X_i\}_{i=1}^{N_1}$ and $\{Y_i\}_{i=1}^{N_2}$ whose log-likelihood functions are denoted by $l_1(X;\mu,\alpha)$ and $l_2(Y;\mu,\beta)$, respectively. The generalized M-estimation problem of interest is defined as

$$\max_{\mu,\alpha,\beta} l_1(X;\mu,\alpha) + l_2(Y;\mu,\beta).$$

A generalization of the subproblem can be cast in terms of the log-likelihood functions directly. Assume there exist two monotone (decreasing) transformations $T_X, T_Y : \mathbb{R} \to \mathbb{R}$ such that $T_X(l_1(X;\cdot,\cdot))$ and $T_Y(l_2(Y;\cdot,\cdot))$ are convex functions. The subproblems, analogous to the EMEP, are

$$h_X(u_1) = \min_{\mu,\alpha,\beta} \{T_Y(l_2(Y;\mu,\beta)) : T_X(l_1(X;\mu,\alpha)) \leq u_1\}.$$

The geometric results in Section 3 still hold with minor modifications. Moreover, under the above convexity assumption, the evaluation of $h_X(u_1)$ can be efficiently performed through standard convex programming techniques. Therefore, the convergence results of Section 4 are still valid.

The above framework encompasses the BFP by taking $T_X(z) = \exp(\frac{2}{N_1} z) - 1$ and $T_Y(z) = \exp(\frac{2}{N_2} z) - 1$.

We now give a simple example of the application of the methodology described above to a Behrens–Fisher-like Problem.

EXAMPLE 6.1. Assume $X \sim N(\mu, \Sigma)$ but, differently from the BFP, $Y$ follows a multivariate Laplacian distribution, that is,

$$f_Y(y) = c_L \exp(-\|y - \mu\|),$$

where $c_L$ is the normalization constant and $\|\cdot\|$ is the Euclidean norm. The related lifted problem can be cast as

$$\begin{aligned}
\min_{\mu,u_1,u_2} \quad & f(u_1, u_2) = \frac{N_1}{2} \log(u_1) + u_2, \\
& u_1 \geq 1 + \mathcal{M}_{\bar{X}}(\mu), \\
& u_2 \geq \sum_{i=1}^{N_2} \|Y_i - \mu\|.
\end{aligned} \quad (21)$$

---

[3]Higher-dimensional convex domains would impose an additional burden in terms of computational complexity.

Here, the problem objective function is concave $(u_1, u_2)$, and therefore the solution must lie on the border of the convex domain of these variables. Such a domain can be written as

$$(22) \quad \mathcal{K} = \left\{ (u_1, u_2) \in \mathbb{R}^2 : \exists \mu \text{ such that } \begin{array}{l} u_1 \geq 1 + \mathcal{M}_{\bar{X}}(\mu) \\ u_2 \geq 1 + \sum_{i=1}^{N_2} \|Y_i - \mu\| \end{array} \right\}.$$

Moreover, the associated subproblems, using $T_X(z) = \exp(\frac{2}{N_1} z) - 1$ and $T_Y(z) = z$, are convex programming problems and have the form

$$h_X(u_1) = \min_\mu \left\{ \sum_{i=1}^{N_2} \|Y_i - \mu\| : \mathcal{M}_{\bar{X}}(\mu) \leq u_1 \right\}$$

and

$$h_Y(u_2) = \min_\mu \left\{ \mathcal{M}_{\bar{X}}(\mu) : \sum_{i=1}^{N_2} \|Y_i - \mu\| \leq u_2 \right\}.$$

Both these problems can be solved via convex quadratic programming, which can be done quite efficiently even in high-dimensional cases.

## APPENDIX A: NOTATION OF CONVEX ANALYSIS

Herein we gather the definitions of relevant concepts in Convex Analysis for this work. We refer to [24] for an analytic exposition of Convex Analysis and to [13] for a more geometric one.

A set $S$ is convex if, for any $x, y \in S$, $\alpha \in [0, 1]$, $\alpha x + (1 - \alpha) y \in S$. An extreme point of a convex set is a point that cannot be written as a strictly $(\alpha < 1)$ convex combination of any other distinct points in the set. A set $P$ is said to be polyhedral if $P = \{x \in \mathbb{R}^n : Ax \leq b\}$, where $A$ is a matrix, and $b$, a vector. It follows that polyhedral sets are convex and their extreme points are its corners. The recession cone $C_S$ of a convex set $S$ is the set of directions that go to infinity in $S$, formally, $C_S = \{d : d + S \subset S\}$.

A function $g : \mathbb{R}^n \to \mathbb{R}$ is said to be convex if, for any $x, y \in \mathbb{R}^n$, and $\alpha \in [0, 1]$, $g(\alpha x + (1 - \alpha) y) \leq \alpha g(x) + (1 - \alpha) g(y)$. A function $f : \mathbb{R}^n \to \mathbb{R}$ is quasi-concave if, for any $x, y \in \mathbb{R}^n$, and $\alpha \in [0, 1]$, $f(\alpha x + (1 - \alpha) y) \geq \min\{f(x), f(y)\}$, or equivalently, the upper level sets of $f$ are convex sets.

Given a convex function $g : \mathbb{R}^n \to \mathbb{R}$, we can define its subdifferential at $x$ as $\partial g(x) = \{s \in \mathbb{R}^n : g(y) \geq g(x) + \langle s, y - x \rangle, \text{ for all } y \in \mathbb{R}^n\}$. The elements of the subdifferential, also called subgradients, play the role of the gradient in case $g$ is nondifferentiable. Note that $\partial g(x)$ is always nonempty.

## APPENDIX B: SOLVING THE EMEP

Consider the convex problem in (7). There are a variety of "general purpose" convergent algorithms that can solve it. Here, we propose a specific algorithm tailored for the particular structure of the EMEP.

Let $\lambda$ be the (nonnegative) Lagrange multiplier associated with the inequality constraint. The first order conditions are necessary and sufficient, and are given by

$$2S_2^{-1}(\bar{Y} - \mu) + 2\lambda S_1^{-1}(\bar{X} - \mu) = 0$$

and

$$\lambda(\mathcal{M}_{\bar{X}}(\mu) - v_1) = 0.$$

Assuming that $\lambda > 0$ (otherwise, the solution is just $\widehat{\mu} = \bar{Y}$), the optimal $\widehat{\mu}$ is a function only of $\lambda$:

(23) $$\widehat{\mu}(\lambda) = (S_2^{-1} + \lambda S_1^{-1})^{-1}(S_2^{-1}\bar{Y} + \lambda S_1^{-1}\bar{X}).$$

Therefore, in order to solve the EMEP, it suffices to compute a root $\lambda^*$ of the nonlinear univariate function

(24) $$m(\lambda) = \mathcal{M}_{\bar{X}}(\widehat{\mu}(\lambda)) - v_1.$$

The algorithm we propose here is based upon the algorithm in Ye [38], who in turn built upon earlier work by Smale [29].

Our algorithm is made up of two main parts. The first part consists of a binary search over intervals of increasing length to find which interval $I_{i^*}$ contains what Smale [29] calls an *approximate root*.

DEFINITION B.1. A point $\lambda^0$ is said to be an approximate root of an analytic real function $m: \mathbb{R} \to \mathbb{R}$ if

$$|\lambda^{k+1} - \lambda^k| \leq (1/2)^{2^{k-1}-1}|\lambda^1 - \lambda^0|.$$

In the second part of the algorithm, Newton's method is used over the interval $I_{i^*}$ to find the approximate root $\lambda^*$. For the sake of exposition, we focus on the case of $m: \mathbb{R} \to \mathbb{R}$ (the results in [38] hold in much greater generality, though). Recall that the Newton iterate for a function $m$ from a current point $\lambda^k$ is

$$\lambda^{k+1} = \lambda^k - \frac{m(\lambda^k)}{m'(\lambda^{k+1})}.$$

Newton's Method (NM) converges quadratically from the very first iteration. In [29], Smale gives sufficient conditions under which a particular point is an approximate root. Although it is hard to verify Smale's condition in general, Ye provided a constructive method to find such a point for a particular class of functions. Ye's results in [38] apply in our case. We now write out Ye's algorithm and prove a complexity result for it in the context of the BFP.

| Binary Search and Newton Method |
| --- |
| **Input:** Upper and lower bounds on the value of the root $[a,b]$, $b \geq a \geq \delta$, tolerance $\delta > 0$. |
| **Step 1.** Define a partition of [a,b] through intervals of the form $I_i = [a(1+1/12)^i, a(1+1/12)^{i+1})$. |
| **Step 2.** Perform binary search on these intervals to find $I_{i^*}$ that contains the true root $\lambda^*$. |
| **Step 3.** Let $\lambda^0 = a(1+1/12)^i$, $k = 0$. |
| **Step 4.** Perform Newton's method from $\lambda^k$: $\lambda^{k+1} \leftarrow \lambda^k - \frac{m(\lambda^k)}{m'(\lambda^k)}$. |
| **Step 5.** Stop if $k > 1 + \log_2(1 + \max\{0, \log_2(b/\delta)\})$ steps. |
| **Step 6.** Else set $k \leftarrow k+1$, and goto Step 4. |

THEOREM B.1.  *After the computation of a spectral decomposition of the matrix $S_1^{1/2} S_2^{-1} S_1^{1/2}$, and given a desired precision $\delta > 0$ and an upper bound $b$ for the solution, the algorithm finds a $\delta$-approximate solution $\widehat{\lambda}$ such that $|\lambda^* - \widehat{\lambda}| < \delta$ in at most*

$$O\left(d \log \log \frac{b}{\delta}\right)$$

*arithmetic operations.*

PROOF. Making the following change of variables/notation

$$w := S_1^{-1/2}(\mu - \overline{X}), \qquad M := S_1^{1/2} S_2^{-1} S_1^{1/2} = PDP^T,$$

$$v = 2 S_2^{-1} S_1^{1/2}(\bar{Y} - \bar{X}) \quad \text{and} \quad s = P^T v,$$

problem (15) is equivalent to

$$h(v_1) = \min w^T M w - v^T w,$$

$$\|w\|^2 \leq v_1$$

up to a constant value (which does not matter for the optimization).

Under the new notation, we can rewrite the function $m$ as

$$m(\lambda) = s^T (D + \lambda I)^{-2} s - v_1 = \sum_{i=1}^{d} \frac{s_i^2}{(D_i + \lambda)^2} - v_1.$$

The function $m(\lambda)$ is analytic and its derivatives can be easily computed as

$$m^{(k)}(\lambda) = (-1)^k (k+1)! \sum_{i=1}^{d} \frac{s_i^2}{(D_i + \lambda)^{k+2}}.$$

Note that $m' < 0$ and $m'' > 0$ (i.e., $m$ is decreasing and convex). Thus, we can evaluate $m$ and $m'$ in $O(d)$ operations. This implies that each Newton step can be implemented in $O(d)$ arithmetic operations.

Let $\lambda^0 = a(1 + 1/12)^{i^*}$ be the left endpoint of the interval selected by binary search. From Ye [38], it follows that $\lambda^0$ satisfies Smale's sufficient condition to be an approximate root. Therefore, NM converges quadratically from the very first iteration (i.e., from $\lambda^0$). From the convexity of $m$, the convergence is monotone, that is, $0 < \lambda^0 < \lambda^k < \lambda^{k+1} < \lambda^* \leq b$ for every $k$ (in particular, we have $|\lambda^1 - \lambda^0| < b$). This implies that we need at most

$$k = 1 + \log_2(1 + \max\{0, \log_2(b/\delta)\})$$

Newton steps to achieve $|\lambda^k - \lambda^*| < \delta$. Moreover, the total number of subintervals is $\frac{1}{\log(1+1/12)} \log(b/a)$. The binary search can thus be implemented in $O(\log \log(b/a))$. The result follows by noting that we can take $a \geq \delta$. $\square$

REMARK B.1. Even when we need to solve the EMEP for many different levels of the Mahalanobis distance function, the spectral decomposition of $S_1^{1/2} S_2^{-1} S_1^{1/2}$ needs to be performed only once. This feature of the algorithm makes it a good auxiliary method for the CLA.

The following lemma illustrates how to obtain subgradients for the function $h_X$ with no additional computational effort, which is of interest for the CLA.

LEMMA B.1. Let $\lambda^*$ be a root of the function $m$ as defined in (24). Then $-\lambda^*$ is a subgradient of $h_X$ at $v_1$.

PROOF. Recall $m(\lambda^*) = 0$ implies that $\mu(\lambda^*)$ minimizes $\mathcal{M}_{\bar{Y}}(\mu) + \lambda^* \mathcal{M}_{\bar{X}}(\mu)$. For any $v$, we have

$$h_X(v_1) = \mathcal{M}_{\bar{Y}}(\widehat{\mu}(\lambda^*)) = \mathcal{M}_{\bar{Y}}(\widehat{\mu}(\lambda^*)) + \lambda(\mathcal{M}_{\bar{X}}(\widehat{\mu}(\lambda^*)) - v_1)$$
$$= \mathcal{M}_{\bar{Y}}(\widehat{\mu}(\lambda^*)) + \lambda^*(\mathcal{M}_{\bar{X}}(\widehat{\mu}(\lambda^*)) - v) + \lambda^*(v - v_1)$$
$$\leq h_X(v) + \lambda^*(v - v_1).$$

Here, we used weak duality ($\min \max \geq \max \min$) as follows:

$$h_X(v) = \min_{\mu} \max_{\lambda \geq 0} \mathcal{M}_{\bar{Y}}(\mu) + \lambda(\mathcal{M}_{\bar{X}}(\mu) - v)$$
$$\geq \max_{\lambda \geq 0} \min_{\mu} \mathcal{M}_{\bar{Y}}(\mu) + \lambda(\mathcal{M}_{\bar{X}}(\mu) - v)$$
$$\geq \mathcal{M}_{\bar{Y}}(\widehat{\mu}(\lambda^*)) + \lambda^*(\mathcal{M}_{\bar{X}}(\widehat{\mu}(\lambda)) - v).$$

Therefore, for every $v$, we have
$$h_X(v_1) - \lambda^*(v - v_1) \leq h_X(v),$$
which implies that $-\lambda^* \in \partial h_X(v_1)$. $\square$

## APPENDIX C: THE DISCRETIZATION ALGORITHM (DA)

Consider the problem (5) for a fixed value of $u_1 = \bar{u}_1$. In this case, the computational problem reduces exactly to solving the EMEP with respect to $X$ at a fixed squared distance level $\bar{u}_1 - 1$. As shown in Section 2, such a problem can be solved directly with the algorithm proposed in Appendix B.

Therefore, given the desired precision, one can discretize the range of the variable $u_1$, $[\bar{L}_1, \bar{U}_1]$, and solve the EMEP for each one of these values. Such a scheme yields the following algorithm.

| Discretization Algorithm |
|---|
| **Input:** Relative tolerance $\varepsilon > 0$, $u_1^1 = (1 + 2\varepsilon/N_1)\bar{L}_1$, $k = 1$. |
| **Step 1.** Evaluate $u_2^k = g(u_1^k)$ and compute $f^k = \frac{N_1}{2}\log(u_1^k) + \frac{N_2}{2}\log(u_2^k)$. |
| **Step 2.** If $(1 + 2\varepsilon/N_1)u_1^k > \bar{U}_1$, compute $f^{k+1} = \bar{U}_1\bar{L}_2$, goto Step 4. |
| **Step 3.** Else set $u_1^{k+1} \leftarrow (1 + 2\varepsilon/N_1)u_1^k$, $k \leftarrow k + 1$, goto Step 1. |
| **Step 4.** Report $\min_{1 \leq i \leq k} f^i$ and the correspondent pair $(\hat{u}_1^{i*}, \hat{u}_2^{i*})$. |

The following complexity result holds for the Discretization Algorithm.

THEOREM C.1. *The Discretization Algorithm reports an $\varepsilon$-solution for the original problem after exactly $\lceil \log(\bar{U}_1/\bar{L}_1)/\log(1 + 2\varepsilon/N_1) \rceil$ loops.*

PROOF. Let $u^* = (u_1^*, u_2^*)$ be a optimal solution. There exists a $k$ such that $u_1^k < u_1^* < (1 + 2\varepsilon/N_1)u_1^k$. We consider $f^{k+1}$ as our candidate. We have

$$f^* = \frac{N_1}{2}\log(u_1^*) + \frac{N_2}{2}\log(u_2^*)$$

(25) $$\leq f^{k+1} = \frac{N_1}{2}\log(1 + 2\varepsilon/N_1) + \frac{N_1}{2}\log(u_1^k) + \frac{N_2}{2}\log(u_2^{k+1})$$

$$\leq \varepsilon + \frac{N_1}{2}\log(u_1^k) + \frac{N_2}{2}\log(u_2^{k+1}) = \varepsilon + f^*,$$

where we also used that $u_2^{k+1} \leq u_2^*$, since $g$ in (14) is decreasing.

The claim on the number of loops follows by noting that we have $u_1^k = \bar{L}_1(1 + 2\varepsilon/N_1)^k \leq \bar{U}_1$ and by taking logs to bound $k$. $\square$

By choosing a sequence $\varepsilon_k \to 0$, we obtain a sequence of $\varepsilon_k$-solutions that converge to the optimal solution of the BFP. One drawback to this method is that it requires solving the EMEP at every point of the discretization. In practice, such requirement may be cumbersome.

# APPENDIX D: MONTE CARLO STUDY OF SIZE

**Acknowledgments.** We are grateful to Victor Chernozhukov, Donald Richards and Vladas Pipiras for thoroughly reading preliminary versions of this paper, and to Pierre Bonami, Oktay Günlük, Jon Lee, Katya Scheinberg and Andreas Wächer for useful discussions. We also thank two anonymous referees and the associated editor for their comments that helped improve the paper.

TABLE 4
*Monte Carlo study of size for the W, LR, LM and B Tests (runs per entry = 10,000)*

| Small size instances | | | Size of the Test $\alpha$ | | | | | | | | | | | |
|---|---|---|---|---|---|---|---|---|---|---|---|---|---|---|
| | | | $\alpha = 0.10$ | | | | $\alpha = 0.05$ | | | | $\alpha = 0.01$ | | | |
| $d$ | $N_1$ | $N_2$ | W | LR | LM | B | W | LR | LM | B | W | LR | LM | B |
| 2 | 10 | 20 | 0.160 | 0.133 | 0.102 | 0.092 | 0.106 | 0.077 | 0.046 | 0.046 | 0.039 | 0.020 | 0.005 | 0.009 |
| 2 | 20 | 40 | 0.138 | 0.122 | 0.106 | 0.103 | 0.081 | 0.067 | 0.050 | 0.051 | 0.023 | 0.013 | 0.007 | 0.009 |
| 2 | 40 | 80 | 0.120 | 0.114 | 0.110 | 0.107 | 0.069 | 0.061 | 0.055 | 0.054 | 0.017 | 0.014 | 0.011 | 0.011 |
| 5 | 25 | 50 | 0.171 | 0.133 | 0.098 | 0.090 | 0.101 | 0.073 | 0.047 | 0.046 | 0.035 | 0.019 | 0.005 | 0.008 |
| 5 | 50 | 100 | 0.124 | 0.110 | 0.094 | 0.090 | 0.068 | 0.055 | 0.041 | 0.041 | 0.017 | 0.011 | 0.007 | 0.008 |
| 5 | 100 | 200 | 0.113 | 0.106 | 0.098 | 0.096 | 0.063 | 0.057 | 0.053 | 0.052 | 0.015 | 0.013 | 0.010 | 0.010 |
| 10 | 50 | 100 | 0.175 | 0.131 | 0.094 | 0.084 | 0.102 | 0.072 | 0.044 | 0.041 | 0.035 | 0.018 | 0.008 | 0.008 |
| 10 | 100 | 200 | 0.137 | 0.118 | 0.099 | 0.093 | 0.074 | 0.062 | 0.047 | 0.047 | 0.019 | 0.012 | 0.009 | 0.009 |
| 10 | 200 | 400 | 0.116 | 0.107 | 0.100 | 0.098 | 0.062 | 0.056 | 0.051 | 0.051 | 0.014 | 0.011 | 0.009 | 0.009 |

| Medium size instances | | | Size of the Test $\alpha$ | | | | | | | | | | | |
|---|---|---|---|---|---|---|---|---|---|---|---|---|---|---|
| | | | $\alpha = 0.10$ | | | | $\alpha = 0.05$ | | | | $\alpha = 0.01$ | | | |
| $d$ | $N_1$ | $N_2$ | W | LR | LM | B | W | LR | LM | B | W | LR | LM | B |
| 25 | 125 | 250 | 0.196 | 0.141 | 0.096 | 0.077 | 0.118 | 0.079 | 0.043 | 0.034 | 0.035 | 0.018 | 0.007 | 0.007 |
| 25 | 250 | 500 | 0.137 | 0.109 | 0.088 | 0.078 | 0.071 | 0.056 | 0.039 | 0.036 | 0.019 | 0.014 | 0.008 | 0.007 |
| 25 | 500 | 1000 | 0.120 | 0.110 | 0.099 | 0.091 | 0.064 | 0.055 | 0.049 | 0.046 | 0.015 | 0.011 | 0.009 | 0.008 |
| 50 | 250 | 500 | 0.232 | 0.158 | 0.096 | 0.065 | 0.144 | 0.089 | 0.040 | 0.027 | 0.041 | 0.018 | 0.005 | 0.005 |
| 50 | 500 | 1000 | 0.147 | 0.117 | 0.091 | 0.079 | 0.083 | 0.061 | 0.044 | 0.038 | 0.021 | 0.013 | 0.008 | 0.007 |
| 50 | 1000 | 2000 | 0.126 | 0.111 | 0.100 | 0.092 | 0.070 | 0.062 | 0.053 | 0.049 | 0.016 | 0.012 | 0.011 | 0.010 |
| 75 | 375 | 750 | 0.262 | 0.170 | 0.098 | 0.064 | 0.167 | 0.098 | 0.048 | 0.029 | 0.059 | 0.025 | 0.008 | 0.004 |
| 75 | 750 | 1500 | 0.166 | 0.131 | 0.097 | 0.084 | 0.098 | 0.072 | 0.048 | 0.038 | 0.025 | 0.016 | 0.009 | 0.007 |
| 75 | 1500 | 3000 | 0.133 | 0.119 | 0.102 | 0.092 | 0.073 | 0.064 | 0.053 | 0.049 | 0.018 | 0.014 | 0.010 | 0.009 |
| 100 | 500 | 1000 | 0.284 | 0.175 | 0.090 | 0.054 | 0.179 | 0.097 | 0.043 | 0.025 | 0.060 | 0.025 | 0.007 | 0.004 |
| 100 | 1000 | 2000 | 0.175 | 0.134 | 0.101 | 0.076 | 0.104 | 0.071 | 0.047 | 0.036 | 0.026 | 0.016 | 0.008 | 0.006 |
| 100 | 2000 | 4000 | 0.139 | 0.117 | 0.099 | 0.087 | 0.073 | 0.061 | 0.050 | 0.042 | 0.017 | 0.013 | 0.009 | 0.007 |

| Large size instances | | | Size of the Test $\alpha$ | | | | | | | | | | | |
|---|---|---|---|---|---|---|---|---|---|---|---|---|---|---|
| | | | $\alpha = 0.10$ | | | | $\alpha = 0.05$ | | | | $\alpha = 0.01$ | | | |
| $d$ | $N_1$ | $N_2$ | W | LR | LM | B | W | LR | LM | B | W | LR | LM | B |
| 200 | 1000 | 2000 | 0.373 | 0.213 | 0.095 | 0.040 | 0.251 | 0.123 | 0.043 | 0.015 | 0.101 | 0.030 | 0.007 | 0.002 |
| 200 | 2000 | 4000 | 0.203 | 0.136 | 0.085 | 0.060 | 0.112 | 0.073 | 0.042 | 0.029 | 0.032 | 0.016 | 0.009 | 0.005 |
| 200 | 4000 | 8000 | 0.153 | 0.128 | 0.099 | 0.085 | 0.084 | 0.064 | 0.049 | 0.039 | 0.019 | 0.014 | 0.010 | 0.007 |

THE FUQUA SCHOOL OF BUSINESS
DUKE UNIVERSITY
1 TOWERVIEW ROAD
PO BOX 90120
DURHAM, NORTH CAROLINA 27708-0120
USA
E-MAIL: abn5@duke.edu

MATHEMATICS DEPARTMENT
TULANE UNIVERSITY
6823 ST. CHARLES AVENUE
NEW ORLEANS, LOUISIANA 70118
USA
E-MAIL: gdidier@tulane.edu